\documentstyle[intlim,amssymb,righttag]{amsart}

\subjclass[2000]{
}


\newtheorem{thm}{\ \ \ Theorem}
\newtheorem{cor}{\ \ \ Corollary}
\newtheorem{lem}{\ \ \ Lemma}

\newtheorem{que}[thm]{Question}
\newtheorem{conj}[thm]{Conjecture}
\theoremstyle{definition}
\newtheorem{defn}{\ \ \ Definition}

\theoremstyle{remark}
\newtheorem{rem}{\ \ \ Remark}
\newtheorem{remo}{\ \ \ Remark on the connection between joinning and Markov operators}

\newtheorem{thank}{\ \ \ Acknowledgment}

\let\texthat\^
\renewcommand\^{\widehat }

\let\standardhash=\#
\renewcommand\#{\raise1pt\hbox{$\scriptstyle\standardhash$}}
\def\tower{
\setlength{\unitlength}{1mm}
\begin{picture}(0,35)
\put (0,0) {\framebox (50,20)} \put (3,-4) {\makebox (3,-4)[br]
{$B_{k}$}} \put (16,-6) {\makebox (16,-6)[b]{$\underbrace {
  \quad  \quad   \quad  \quad \quad  \quad  \quad  \quad  \quad  \
 \quad  \quad \quad  \quad
 }_{B_{k-1}}$}}
\put (17,-16) {\makebox (17,-16)[b]{$\overline {
     \quad  \quad  \quad \quad ~ ~ p_{k-1}
\quad  \quad   \quad  \quad \quad   \quad ~ ~
 }$}}
\multiput (5,0)(5,0){10}{\line (0,20){20}} \put (0,4) {\line
(4,0){50}} \put (3,0){\vector (0,1){4}} \thicklines \put (0,0)
{\line (5,0){5}} \multiput (0,21)(0,1){5}{\line (5,0){5}} \put
(-3,27){\makebox (-2,27)[bl] {$a_1^{(k-1)}$}} \put
(5,21){\line(5,0){5}} \put (5,22){\line(5,0){5}} \put
(8,24){\makebox (8,24)[bl] {$a_2^{(k-1)}$}} \put (10,21){\makebox
(10,21)[b] {$\cdots \cdots$}} \multiput (25,21)(0,1){6}{\line
(25,0){5}} \put (23,28){\makebox (23,28)[bl] {$a_j^{(k-1)}$}} \put
(26,21){\makebox (26,21)[b] {$\cdots \cdots$}} \multiput
(45,21)(0,1){3}{\line (45,0){5}} \put (45,25){\makebox (45,25)[bl]
{$a_{p_{k-1}}^ {(k-1)}$}} \put (-40,-25){\centerline{Figure 1 :
$k^{\hbox {th}}$--tower.}}
\end{picture}}

\newcommand\sym{\fam\comfam\com}
\font\tensym=msbm10 at 12pt \font\sevensym=msbm7
\font\fivesym=msbm5 
\newfam\symfam
\textfont\symfam=\tensym \scriptfont\symfam=\sevensym
\scriptscriptfont\symfam=\fivesym
\renewcommand\sym{\fam\symfam\relax}

\newcommand\N{{\sym N}}
\newcommand\Z{{\sym Z}}

\newcommand\T{{\sym T}}

\newcommand\1{{\sym I}}

\title
{On the spectrum of $\alpha$-rigid maps}
\author{E.~H.~EL ABDALAOUI}

\begin{document}


\begin{abstract}
It is shown that there exists an $\alpha$-rigid transformation
with $\alpha$ less or equal to $\frac12$ whose spectrum has
Lebesgue component. 
This answers the question raised by Klemes and
Reinhold in \cite{Klemes-Reinhold}.
Moreover, we introduce a new criterium to identify a large class of $\alpha$-rigid transformations with singular spectrum.
\end{abstract}

\maketitle

\section{Introduction}

In his 1980's paper \cite{Katok2}, A. Katok showed that interval
exchange maps are not mixing. As a consequence of Katok's proof interval
exchange maps are $\alpha$-rigid. In a later work, Veech in \cite{Veech}
investigated the spectral properties of interval exchange maps and
showed that almost every minimal interval exchange maps has
singular simple spectrum. Before, Oseledets in \cite{Oseledets}
proved that for any interval exchange transformation the maximal
spectral multiplicity is bounded above by $r-1$, where $r$ is the
number of intervals exchanged. Moreover, he constructed the first
example of a transformation with continuous spectrum and finite
multiplicity greater than $1$. Since the example is an exchange of
$30$ intervals, the maximal spectral multiplicity $m$ satisfies $2
\leq m \leq 29$. Robinson \cite{Robinson} constructed ergodic
interval exchange transformations with arbitrary finite maximal
spectral multiplicity. It follows from \cite{Katok3} that all the
examples constructed by Oseledets and Robinson have singular
spectrum. However Katok's theorem \cite{Katok2} remains at the present time the
only universal result about the spectrum of interval exchange maps.
Following Veech one may ask if there is some other ``universal"
spectral property satisfied by interval exchange maps. More
precisely, as in \cite{elabdal} one may ask the following question.

\begin{que}
Does any interval exchange map have singular spectral type?
\end{que}

\noindent{}The answer to the above question is affirmative in the
case of three interval exchange maps. Note that using the result in
\cite{guenais} one may show that every ergodic interval exchange
transformation on three intervals has singular spectrum.\\

\noindent{}On the other hand, Klemes and Reinhold in
\cite{Klemes-Reinhold} proved that for any $\alpha \in ]0,1[$,
there exists an $\alpha$-rigid rank one transformation with
singular spectrum and asked the following question.

\begin{que}
Does any $\alpha$-rigid transformation have singular spectral
type?
\end{que}

\noindent{}It is well known that the answer is affirmative if $\alpha >\frac12$ (Baxter's theorem).
Note that if this was true for all $\alpha \leq 1$  then the answer the question
$1$ would be affirmative too.

In this paper we shall prove that this is not the case. More precisely, we prove in section 3 that the
Mathew-Nadkarni transformation (which has Lebesgue component in the spectrum) is $\displaystyle\frac12$-rigid. We
recall that Mathew-Nadkarni introduced this transformation in 1984
in \cite{Mathew-Nadkarni} to answer the question raised by Helson
and Parry \cite{Helson-Parry} of whether there exists an ergodic
measure preserving transformation which has a Lebesgue component
in its spectrum with finite non-zero multiplicity. Helson and
Parry also mentioned the problem, attributed to Banach, whether
there exists an ergodic measure preserving transformation on a
finite measure space whose spectrum is simple Lebesgue. In
\cite{Rokhlin}, Rokhlin mentioned the problem of finding an
ergodic measure preserving transformation on a finite measure
space whose spectrum is of Lebesgue type with finite multiplicity.
Another contribution to the Banach-Rohklin question is due to M.
Queff\'elec \cite{Queffelec1} who proved that the spectrum of the
Rudin-Shapiro substitution has Lebesgue component of multiplicity
two. It turns out that the Rudin-Shapiro substitution is
$\displaystyle\frac12$-rigid.

In section 4 we will present other examples of
dynamical systems arising from substitutions with Lebesgue
component of multiplicity two whose rigidity constant is less than $\displaystyle\frac12$.

It is an easy exercise to prove that $\alpha$-rigid
transformations are not mixing. Since $\alpha$-rigidity does not
imply that the spectrum is singular (as we shall see in the section 3), the absence
of mixing does not imply that the spectrum is singular. Let us
mention that any aperiodic measure-preserving transformation can
be realized as an exchange of an infinite number of intervals \cite{Arnoux}.\\

Our work on the question of Klemes and Reinhold (question 2) includes a survey
of the results of Dekking-Keane \cite{Keane} and Queff\'elec
\cite{Queffelec2}. F. M. Dekking and M. Keane showed that the
dynamical systems arising from substitutions are not mixing. The
ingredients of the proof lead to establish that the dynamical systems
arising from substitutions are $\alpha$-rigid. The procedure to
check the constant $\alpha$ will be presented in the last section
without proof. Queff\'elec showed that the substitution which
gives rise to the Rudin-Shapiro sequence has Lebesgue component of
multiplicity two together with a discrete component. (Kamae
\cite{Kamae} had earlier shown that the
correlation measure of the Rudin-Shapiro sequence is Lebesgue).\\

In this section 2 we shall exhibit a large class of
$\alpha$-rigid transformations with singular spectrum. More
precisely, we shall prove that the transformation satisfying the
restricted Beurling condition is singular. We say that the
transformation $T$ satisfies the {\bf {restricted Beurling
condition}} if the following holds
\begin{eqnarray*}
\displaystyle &&\left\lbrace \ \displaystyle \sum_{i \in
\Z}a_iU_T^i ~~:~~ a_i > 0
  {\rm {~~for~~some~}} i {\rm {~~and~~}}
  \displaystyle \sum_{n \geq 0}\frac{\log(\displaystyle \sum_{k \leq -n}a_k^2)}{n^2}=-\infty\right\rbrace
  \\ &&\bigcap
   \left ( {\overline{\left\lbrace U_T^n, n\in \Z\right\rbrace}}^{W}
  \setminus \left ({\left\lbrace U_T^n, n\in \Z\right\rbrace} \right ) \right)
  \neq
  \emptyset.
\end{eqnarray*}
 \noindent{}Where $U_T$ is the operator defined by $U_T
f(x)=f(T^{-1}x)$ and ${\overline{\left\lbrace U_T^n, n\in \Z\right\rbrace}}^{W}$
 is the weak closure power of $T$.\\

We mention that the $\alpha$-rigid rank one
transformations constructed by Klemes \& Reinhold statify the condition above. Thus, our result can be considered as a generalization
of the Klemes-Reinhold result. There is another motivation for our
work. In 1979 D. Rudolph \cite{Rudolph} introduced the notion of
minimal self-joinings as the foundation for a machinery that
yields a wide variety of counterexamples in ergodic theory. But
the property of minimal self-joinings is meager with respect to
the weak topology on the group of the all automorphisms. Later, in
1983 A. del Junco and M. Lema\'{n}czyk in \cite{Del} showed that
these constructions, as well as many others in \cite{Rudolph}, can
be based on a much weaker property that is in fact generic
(residual in the weak topology). A special case of this property is the property of $\kappa$-mixing (see Katok \cite{Katok1}, Stepin \cite{Stepin}) which implies
the mutual singularity of the convolution powers of the maximal
spectral type of an automorphism $T$. We remark also that any $\kappa$-mixing transformation is
$(1-\kappa)$-rigid.

\begin{que}\label{powers}
Are the convolution powers of the maximal spectral type of any
non-mixing transformation with minimal self-joinings property
pairwise singular?
\end{que}

This would imply that the spectrum of any non-mixing transformation with minimal self-joinings is singular. In this paper we are not able to answer the question \ref{powers}. Nevertheless we shall exhibit in section 2 a large class of non-mixing transformation with minimal self-joinings with singular spectrum.

Recently, some progress was obtained by Prikhod'ko and Ryzhikov in
\cite{Prikhod'ko-Ryzhikov}. The authors show that the well-known
Chacon transformation \cite{Chacon2} possesses this property. It is well known
that the Chacon transformation has the minimal self-joining
property \cite{Deljunco}. This result can be extended easily to
the case of staircase transformations with bounded cutting
parameter. The latter examples include as a special case the
Klemes-Reinhold examples of
$\alpha$-rigid rank one transformations. \\

We recall now some basic facts on
spectral theory. A nice account can be found in the appendix of \cite{Parry}.

\subsection*{1.1. Spectral measures}
Given a measure preserving invertible
transformation $T:X\mapsto X$  and denoting as above by $U_T f$ the operator  $U_T
f(x)=f(T^{-1}x)$, for any $f\in L^2(X)$ there
exists a positive measure $\sigma_f$ on the unit circle $S^1$ defined by $\hat{\sigma}_f(n)= <U_T^nf,f>$. \\

\begin{defn}
The maximal spectral type of $T$ is the equivalence class of Borel
measures $\sigma$ on $\T$ (under the equivalence relation $\mu_1
\sim \mu_2$ if  and only if $\mu_1<<\mu_2$ and $\mu_2<<\mu_1$),
such that
 $\sigma_f<<\sigma$ for all $f\in L^2(X)$ and
if $\nu$ is another measure for which $\sigma_f<<\nu$
for all $f\in L^2(X)$ then $\sigma << \nu$.\\
\end{defn}

By the canonical decomposition of $L^2(X)$ into decreasing cycles
with respect to the operator $U_T$, there exists a Borel measure
$\sigma=\sigma_f$ for some $f\in L^2(X)$, such that $\sigma$ is in
the equivalence class defining the maximal spectral type of $T$.
By abuse of notation, we will call this measure the maximal
spectral type measure. It can be replaced by any other measure in
its equivalence class. The reduced maximal type $\sigma_{T}^{(0)}$ is the
maximal spectral type of $U_T$ on $L_0^2(X)\stackrel{\rm
{def}}{=}\{f \in L^2(X)~:~ \displaystyle \int f d\mu=0 \}$. The
spectrum of $T$ is said to be discrete (resp. continuous, resp.
singular, resp. absolutely continuous , resp. Lebesgue ) if
$\sigma_{T}^{(0)}$ is discrete ( resp. continuous, resp. singular, resp.
absolutely continuous with respect to the Lebesgue measure or
equivalent to the Lebesgue measure).  We write
\[Z(h) \stackrel{\rm {def}}{=} \overline {{\rm {span}} \{U_T^nh,n \in \Z \} }.
\]
The transformation $T$ is said to have simple spectrum, if there exists $h \in L^2(X)$
such that
\[Z(h) =L^2(X).\]
Two dynamical systems $(X,{\mathcal{A}},\mu,T)$ and $(Y,{\mathcal{B}},\nu,S)$ are spectrally disjoint if $\sigma_{T}^{(0)}$ and $\sigma_{S}^{(0)}$ are
mutually singular.

\subsection*{1.2. $\alpha$-rigid transformations}

\begin{defn}
A measure-preserving transformation $T$ on the probability space
$(X, {\mathcal {B}}, \mu)$ is said to be $\alpha$-rigid, where
$\alpha \in ]0,1]$, if there is a sequence of integers
$\{n_k\}_{k \in \N}$ such that
$$ \lim_{k \rightarrow\infty}\mu(T^{n_k}A \cap A) \geq \alpha \mu(A),
~~~\forall A \in {\mathcal{B}}.$$
\end{defn}

The notion of $\alpha$-rigidity has been formulated in 1987 by N. Friedman
\cite{Friedman1}. Besides this, in 1969 J. Baxter proved
\cite{Baxter}

\begin{thm} [Baxter]
The spectrum of any $\alpha$-rigid transformation is singular if
$\alpha > \displaystyle \frac12$.
\end{thm}

Note that this theorem is a consequence of the following result

\begin{thm}[Ryzhikov\cite{Ryzhikovemail}]
Any $\alpha$-rigid transformation $T$ with
$\alpha > \displaystyle \frac12$ is spectrally disjoint from any mixing transformation $S$
\end{thm}

\begin{pf}
By assumption, there exists a sequence of integers and a Markov operator $P$ such that $(U_T^{n_i})$ converges weakly to $\alpha I +(1-\alpha) P$. Let $J$ be any operator such that
$$U_T J= J U_S.$$
\noindent{}Then, for any $f \in L_0^2(X)$ we have
$$\alpha Jf +(1-\alpha) PJf=0.$$
Since  the norm of $P$ is $\leq 1$, the operator $\alpha I +(1-\alpha) P$ is invertible. Hence,  $Jf$ must be $0$ thus $J\equiv 0$. It follows that $T$ is spectrally disjoint from $S$ and the proof of the theorem is complete.
\end{pf}

\begin{rem}
The proof above gives more, namely if the weak closure of the powers of the operator $U_T$ contains any invertible operator, then $T$ is spectrally disjoint from any mixing maps. In fact,
assume that there exists an invertible operator $V$ in the weak closure of the power of $T$ which means that there exists a sequence of integers $(n_i)$ such that
$$ \lim_{i \longrightarrow +\infty}(U_T^{n_i}f,f)=(Vf,f),$$
\noindent{}for any $f \in L^2_0(X)$. Let $S$ be any mixing map. By the Lebesgue decomposition of
$\sigma_T^{(0)}$ with respect to the maximal spectral type of $S$ we have
$$ \sigma_T^{(0)}= \sigma_s+\sigma_a.$$
It is well known that there exists a function $f \in L^2_0(X)$ such that $\sigma_f=\sigma_a$. But for any function $g \in L^2_0(X)$ we have
$$\lim_{j\longrightarrow +\infty} \lim_{i\longrightarrow +\infty}(U_T^{n_i-n_j}g,g)=(Vg,Vg).$$
Applying the Riemann-Lebesgue theorem we get
$$(Vf,Vf)=0,$$
therefore $f=0$ and the proof is complete.
\end{rem}

\subsection*{1.3. Joinings}

Let $T$ be any ergodic automorphism of a Lebesgue space
$(X,{\mathcal B},\mu)$. The  centralizer $C(T)$ of $T$ is the
semi-group of all endomorphisms $S~ :~ (X,{\mathcal B},\mu)
\rightarrow (X,{\mathcal B},\mu)$ such that $ST=TS$. We say that $\lambda$ is a $n$-joining of $T$  if  $\lambda$  is a $T^{\otimes n}$-invariant probability~measure on ${\mathcal {B}}^{\otimes n}$
and  $\lambda_{|{\mathcal {B}}}=\mu,$
where $T^ {\otimes n}\stackrel{\rm {def}}{=} T \times T \cdots \times T$. Denote by
$J(T,\cdots,T)$ the space of all $n$-joinings of $T$.\\

A standard example of ergodic 2-joinings comes from the
centralizer of $T$. More precisely, if $S \in C(T)$ then the measure
given by
\begin{eqnarray*}
\Delta_S(A \times B)=\mu(A \cap S^{-1}B)
\end{eqnarray*}
is a 2-joining. For $S=T^n$, $n \in \Z$, we put
$\Delta_n \stackrel {\rm {def}}{=}\Delta_{T^n}$ and $\Delta \stackrel {\rm {def}}{=}\Delta_0.$\\

Following \cite{Del-Ru} and \cite{Veech2}, $T$ is called $2-${\it
fold} simple if each ergodic 2-joining of $T$ is either on the
graph of some $S \in C(T)$ or is the product
measure $\mu \times \mu$. It is easy to see that 2-fold
simplicity of $T$ implies that $C(T)$ is a group (consider $\mu_S(A \times B)=\mu(S^{-1}A \cap B)$).\\

Now, for $n>1$ and any $S_i \in C(T)$, $i=1,\cdots,n,$ the $n$-joining measure $\mu_{S_1,\cdots,S_n}$ given by
\[
\mu_{S_1,\cdots,S_n}(A_1\times \cdots \times A_n)=\mu(S_1^{-1}A_1
\cap \cdots \cap S_n^{-1}A_n),\]
is said to be {\it off-diagonal}.\\

The transformation $T$ is said to be {\it simple} if $C(T)$ is a group and for
every $n \geq 2$ and every ergodic $n-$joining $\lambda$ the set
$\{1,\cdots,n\}$ can be split into subsets $s_1,\cdots,s_k$ such
that each ${\lambda|}_{{\mathcal {B}}^{s_i}}$ is off-diagonal and
$\lambda$ is their product \cite{Del-Ru}. If $T$ is
2-fold simple and $C(T)$ is trivial, i.e., $C(T)=\{T^i ~:~i
\in\Z\}$, then $T$ is said to have the property of minimal
self-joining of order 2 ($T \in MSJ(2))$. Moreover $T$ is said to have
minimal self-joinings of any order if $T$ is simple and $C(T)$ is
trivial. The class of $MSJ(2)$ transformations is contained in the
class of transformations with finite joining rank. The joining
rank of any ergodic transformation $T$, written $\rm {jrk}(T)$, is
defined as the minimum of $r \in \N$ for which each ergodic
$r-$fold self-joining $\nu$ has some pair $i<j$ such that the
two-dimensional marginal $ \nu_{ {\mathcal {B}}^{(i)} \times
{\mathcal {B}}^{(j)} }$ is trivial (i.e.,
 $ \nu_{{\mathcal {B}}^{(i)} \times {\mathcal {B}}^{(j)} }  \in
\{\mu \times \mu\} \cup \{\Delta_n\}_{n \in \Z}$). The standard
examples of transformations with finite joining rank greater than 2
is given by the powers of any weak mixing transformation $T$ in
$MSJ(2)$. In addition, a non-mixing transformation with the minimal
self-joining property is $\alpha-$rigid. Indeed,
we have the following
\begin{lem}[\cite {Ryzhikov}]
Let $T$ be a non-mixing map in $MSJ(2)$. Then $T$ is $\alpha$-rigid
for some $\alpha \in ]0,1[$.
\end{lem}

\begin{pf}
Since the set $J(T,T)$ is compact and $T$ is not mixing, there is
a sequence $n_k\rightarrow\infty$ such that the sequence of
diagonal measures $\Delta_{n_k}$ converges to some measure
$\lambda \in J(T,T)$, where $\lambda \neq \mu \times \mu$. Since
$T \in MSJ(2)$, by the ergodic decomposition, we have
\begin{eqnarray}
\lambda=\beta\mu \times \mu+\sum_{i \in\Z}a_i\Delta_i,
\end{eqnarray}
where $a_i \geq 0$ and $\sum_{i \in
\Z}a_i=1-\beta$. Hence, for some $i$ we have
$$\lambda \geq a_i \Delta_i,~~a_i >0.$$
Thus, for the sequence $m_k=n_k-i$, we obtain
\begin{eqnarray*}
\lim_{k \rightarrow \infty}\mu(T^{m_k}A \cap A) \geq \alpha \mu(A)
\end{eqnarray*}
where $\alpha\stackrel{\rm {def}}{=}a_i$ and the proof of the lemma is complete.
\end{pf}

\noindent{}\begin{remo}
It is well known and easy to see that for any joining $\lambda$ between
two dynamical systems $(X,{\mathcal{A}},T,\mu)$ and $(Y,{\mathcal{B}}, S,\nu)$ there exists a Markov operator
$J~:~ L_0^2(X,\mu) \longrightarrow L_0^2(Y,\nu)$ such that
\begin{enumerate}
  \item $J \circ U_T = U_S \circ J,$
  \item $J\chi_X = \chi_Y$ and $J^*\chi_Y = \chi_X,$
  \item $f > 0$ and $g > 0$ implies $Jf > 0$ and $J^*g > 0,$
  \item for any $f \in L_0^2(X,\mu)$ and $g \in L_0^2(Y,\nu)$ we have
$$\int_{X \times Y}
f(x)g(y)d\lambda(x,y)=
\int_Y Jf(y)g(y)d\nu(y).$$
\end{enumerate}
\noindent Here $J^*$ is the adjoint of $J$ defined by $(Jf,g) = (f, J^*g)$ for $f \in L_0^2(X,\mu)$ and $g \in L_0^2(Y,\nu)$. The correspondence is one-to-one. The systems $(X,{\mathcal{A}},T,\mu)$ and $(Y,{\mathcal{B}}, S,\nu)$ are disjoint (in the sens of Furstenberg) if J is trivial, where trivial means
$Jf = \displaystyle\int fd\mu$, for any $f \in L^2(X,\mu)$. To tackle some problems is ergodic theory using this approach is proposed for example in \cite{RyzhikovS}. Using the same approach, Fra\c czek-Lema\'nczyk \cite{frac-lem1} show that the ergodic $\alpha$-rigid maps with $\alpha > 0$ is disjoint from all mixing maps. Indeed, let $\lambda$
be any ergodic joining between a ergodic $\alpha$-rigid map $T$ and a mixing map $S$ and let $J$ be a associated
Markov operator. By the $\alpha$-rigidity of $T$ it is easy to see that there is a sequence
of integers $(n_k)_{k\in \N}$ and a Markov operator $P$ such that $(U_T^{n_k})$ converges weakly to $\alpha Id+(1-\alpha)P$. Thus
$J U_T^{n_k}$ converges weakly to $\alpha J +(1-\alpha)JP$. Since $J U_T^{n_k}= U_S^{n_k}J$ and $S$ is mixing we have $(U_S^{n_k}Jf,g) = 0$ for each $f \in L_0^2(X,\mu)$ and $g \in L_0^2(Y,\nu)$. As a consequence  $$\alpha(Jf,g) + (1 -\alpha)(JPf,g) = 0$$ for each $f \in L_0^2(X,\mu)$ and $g \in L_0^2(Y,\nu)$. This gives, for all $f \in L^2(X,\mu)$, $$\alpha Jf + (1 -\alpha)JPf=\int f d\mu.$$
\noindent Since the dynamical system $(X \times X,{\mathcal{A}} \bigotimes {\mathcal{B}},T \times S,\mu \times \nu)$ is ergodic and $JPU_T=U_SJP$, it follows that the operator $f \in L^2(X,\mu) \longmapsto \displaystyle \int f d\mu$ is indecomposable (i.e., is an extreme point of the subspace of Markov operators). Therefore $\displaystyle Jf=\int f d\mu$ for each $f \in L^2(X,\mu)$ and the proof is complete.\\
\end{remo}

In the following we shall establish the connection between the old
problem of Banach on the existence of dynamical systems with simple
Lebesgue spectrum and the $MSJ$ property. For that, we recall the
following conjecture.

\begin{conj}\label{parry}(\cite{Parry},pp.~50)
If a reduced maximal spectral type of some transformation with
simple spectrum is absolutely continuous then
it is \linebreak Lebesgue.\\
\end{conj}

\begin{thm}
If $~T~$  is a transformation with finite joining rank and simple
Lebesgue spectrum then Conjecture \ref{parry} implies that $T$ is in the class $MSJ$.
\end{thm}

\begin{pf}
By the King-Thouvenot theorem \cite{King-Thouvenot}, $T$ is the
$e$-extension ($e \in \N^*$) of the power (say $p$) of some transformation $S$
with the $MSJ$ property. But $S^p$ is a factor of $T$. It follows that
$S^p$ and $S$ have simple Lebesgue spectrum. Since
the multiplicity of $S^p$ is $p$ we have $p=1$. Thus $T$ is
an $e$-extension of $S$. Now, by the classical spectral decomposition
of the group extension we get $e=1$ and this concludes the proof of the theorem.
\end{pf}

\section{Main result: Beurling condition.}

In this section we introduce the following condition.
\begin{defn} We say that the transformation $T$ satisfies the ``Beurling condition " if
\begin{eqnarray*}
\displaystyle
&&\left\lbrace \ \displaystyle \sum_{i \in \Z}a_iU_T^i ~~:~~ a_i \neq 0
  {\rm {~~for~~some~}} i {\rm {~~and~~}}
  \displaystyle \sum_{n \geq 0}\frac{\log(\displaystyle \sum_{k \leq -n}a_k^2)}{n^2}=-\infty\right\rbrace
  \\ &&\bigcap
   \left ( {\overline{\left\lbrace U_T^n, n\in \Z\right\rbrace}}^{W}
  \setminus \left ({\left\lbrace U_T^n, n\in \Z\right\rbrace} \right ) \right)
  \neq
  \emptyset.
\end{eqnarray*}
\end{defn}

This definition is inspired by early work of Beurling on a quasi-analytic class of mappings and on the uncertainty principle in harmonic
analysis \cite{Beurling} combined by ideas coming from the joining theory.

The main result of this paper is the following theorem.

\begin{thm}\label{beurling}
If the transformation $T$ satisfy the Beurling condition
then the spectrum of $T$ is singular.
\end{thm}

To prove this theorem we shall need the following lemma. We recall that a measure $\nu$ on the circle is called a {\em
Rajchman measure} if $\displaystyle \lim_{n \to \infty}
\widehat{\nu}(n)=0$.\\

\begin{lem}[Translation lemma]
Suppose $\mu$ is a probability measure on $\T$ and $(n_k)_{k \in
\N}$ is a sequence of distinct integers. Define $\mu_k$ by
\[
\mu_k=e^{i n_k \theta} d\mu(\theta),~ k=1,2,3,\cdots
\]
If $(\mu_k)_{k \in \Z}$ converges in the weak* topology to
$\sigma$, then $\sigma$ is singular with respect to any Rajchman measure $\nu$. In fact, if $\mu=\mu_s+\mu_a$
is the Lebesgue decomposition of $\mu$ with respect to the Rajchman measure $\nu$ on $\T$, then
\[
|\sigma|(E)\leq\mu_s(E), {\rm {~for~every~Borel~set~}} E {\rm
{~in~}} \T.
\]
\end{lem}
\bigskip
The proof of this lemma is the same as the proof given in {\cite[pp.~66]{Rudin}} in the case of Lebesgue measure. We will also need the following Beurling theorem proved in \cite{Beurling}

\begin{thm}[Beurling's Theorem]\label{Tbeurling}
Let $f$ be the function on torus $\T$ define by
\[
f(\theta)=\sum_{j \in \Z}a_j e^{ i j\theta}.
\]
\noindent and assume that the following condition holds
\begin{eqnarray*}
 \sum_{n \geq 0}\frac{\log(\sum_{k \leq -n}a_k^2)}{n^2}=-
\infty.
\end{eqnarray*}
Then if $f=0$ on the set of positive measure we have
$f=0$ a.e. on $\T$.
\end{thm}
\bigskip
\begin{pf}(of Theorem \ref{beurling})
Let $\sigma\stackrel{\rm {def}}{=}\sigma_h$ be the reduced maximal
spectral type of $T$. It follows from the Beurling condition
that there exists a sequence of integers $(n_k)_{k \in \N}$ such
that
\begin{eqnarray}
e^{i n_k \theta}d\sigma \stackrel {\rm {w*}}{\rightarrow}\sum_{j
\in \Z}a_j e^{ij\theta}d\sigma.
\end{eqnarray}

\noindent{}Let $\sigma=\sigma_s+\sigma_a$
 be the Lebesgue decomposition of
$\sigma$ with respect to the Lebesgue measure $m$ on $\T$. Put
 $$f(\theta)\stackrel{\rm {def}}{=}
\sum_{j \in \Z}a_je^{ij\theta}.$$ Then, by the translation lemma, we
have
\begin{eqnarray}
\sigma_a\{x~:~f(x)=0\}=\sigma_a(\T).
\end{eqnarray}
\noindent{}Now Beurling's theorem (Theorem \ref{Tbeurling}) says that
\[
f\neq 0 \Rightarrow m\{x~:~f(x)=0\}=0.
\]
Hence $\sigma_a=0$ and the proof of the
theorem is completes.
\end{pf}

As an immediate consequence of Theorem \ref{beurling} we obtain the following

\begin{cor}
If the transformation $T$ satisfies the following restricted Beurling condition
\begin{eqnarray*}
\displaystyle
&&\left\lbrace \ \displaystyle \sum_{i \in \Z}a_iU_T^i ~~:~~ a_i > 0
  {\rm {~~for~~some~}} i {\rm {~~and~~}}
  \displaystyle \sum_{n \geq 0}\frac{\log(\displaystyle \sum_{k \leq -n}a_k^2)}{n^2}=-\infty\right\rbrace
  \\ &&\bigcap
   \left ( {\overline{\left\lbrace U_T^n, n\in \Z\right\rbrace}}^{W}
  \setminus \left ({\left\lbrace U_T^n, n\in \Z\right\rbrace} \right ) \right)
  \neq
  \emptyset.
\end{eqnarray*}
\noindent{}Then $T$ is $\alpha$-rigid with singular spectrum.
\end{cor}

\begin{rem}\label{disjointmixing}
Note that we have actually proved that $T$ is spectrally disjoint from any mixing maps $S$ provided that
\begin{eqnarray}\label{Ryzmixing}
\sigma_S^{(0)}\left(\left\{x~:~f(x)=0\right\}\right)=0.
\end{eqnarray}
\noindent{}Observe that (\ref{Ryzmixing}) holds if $f$ is an analytic function on the circle.
\end{rem}
\subsection*{2.1. Applications.}

Using Theorem \ref{beurling} we shall give in this section a simple proof of
some well known results on the singularity of the spectrum of some
special maps of rank one maps. More precisely, we shall give a
simple proof of the singularity of Chacon maps and the staircase
maps with bounded cutting parameter. Recall that for any
$\varepsilon >0$ one may construct a staircase maps with bounded
cutting parameter, say $p$, such that the $\alpha$-rigid constant is
smaller that $\varepsilon$. In fact, the
$\alpha$-rigid constant is $\displaystyle \frac1{p}$.\\

\noindent{}Let us recall the definition of rank one maps and
in particular the staircase maps. Using the cutting and stacking method
described in \cite{Friedman2}, \cite{Friedman3}, one can define inductively a family of
measure preserving transformations, called rank one
transformations, as follows. \vskip 0.1cm Let $B_0$ be the unit
interval equipped with the Lebesgue measure. At stage one we
divide $B_0$ into $p_0$ equal parts, add spacers and form a stack
of height $h_{1}$ in the usual way. At the $k^{th}$ stage we
divide the stack obtained at the $(k-1)^{th}$ stage into $p_{k-1}$
equal columns, add spacers and obtain a new stack of height
$h_{k}$. If during the $k^{th}$ stage of our construction  the
number of spacers put above the $j^{th}$ column of the
$(k-1)^{th}$ stack is $a^{(k-1)}_{j}$, $ 0 \leq a^{(k-1)}_{j} \leq
\infty$,  $1\leq j \leq p_{k}$, then we have

$$h_{k} = p_{k-1}h_{k-1} +  \sum_{j=1}^{p_{k-1}}a_{j}^{(k-1)}.$$

\vskip .5 cm \hskip 3.5cm
 \tower
\vskip 3.0cm
\noindent{}Proceeding in this way, we get a rank one transformation
$T$ on a certain measure space $(X,{\mathcal B} ,\nu)$ which may
be finite or
$\sigma-$finite depending on the number of spacers we added. \\

\noindent{} The construction of any rank one transformation thus
needs two parameters: $(p_k)_{k=0}^\infty$ (cutting
and stacking parameter), and $((a_j^{(k)})_{j=1}^{p_k})_{k=0}^\infty$
(spacers parameter). Let

$$T \stackrel {def}= T_{(p_k, (a_j^{(k)})_{j=1}^{p_k})_{k=0}^\infty}.$$
\noindent In the case of staircase maps the spacers parameter is given by
$$a_j^{(k)}=j-1,~~{\rm {for}}~~j=1 \hdots p_k-1~~{\rm
{and}}~~a_{p_k}^{(k)}=0.$$

\noindent{}The classical Chacon map \cite{Chacon2} corresponds to
the case $p_k=3$, for every $k \in \N^*$. It is easy to see that the
Chacon map is $\displaystyle \frac13$-rigid. More generally, it is easy to prove
that the staircase with bound cutting parameter (say $p$) is
$\displaystyle \frac1{p}$-rigid (in fact, for any measurable set $A$ we have
$\liminf \mu(T^{h_n}A \cap A) \geq
\displaystyle \frac1{p}\mu(A)$). The following lemma can be may proved using Theorem \ref{beurling} and the Remark \ref{disjointmixing}. Nevertheless  for the convenience of the reader we present a different proof.\\

\begin{lem}\label{beurlingmixing} Let $T$ be a weak mixing transformation and assume that there is a sequence of integers $(n_k)$ such that $(U_T^{n_k})_{k \in \N}$ converges weakly
to $P(T)$, where $P(z)$ is a nonzero analytic function on the circle. Then $T$ is spectrally disjoint from any mixing transformation $S$.
\end{lem}
\begin{pf}(Ryzhikov \cite{Ryzhikovemail},\cite{frac-lem1})
Let $J$ be any operator such that
$$U_T J= J U_S.$$
\noindent{}Then $P(T)Jf=0$ for each $f \in L^2_0(X)$. Since the maximal spectral type of $T$ is continuous and $P(z)$ is analytic we have
$$\ker{P(T)}=\{0\}.$$
\noindent{}Thus $Jf=0$ for each $f \in L^2_0(X)$ i.e., $J$ is trivial and the proof of the lemma is complete.
\end{pf}

\begin{cor}
The staircase maps with bounded cutting parameter
are spectrally disjoint from any mixing maps.
\end{cor}
\begin{pf}
By the definition of staircase maps it is easy to prove that the
sequence $T^{h_n}$ converges weakly to $\displaystyle{\
\sum_{j=0}^{r-1}\frac1{r-1}T^{j}}$ where $T^0$ is the identity map.
\end{pf}

\noindent{}One may apply Lemma \ref{beurlingmixing} also to the ``historical''
Chacon map \cite{Chacon1}, given by
$$p_k=2,~~~~k \in \N^* {\rm {~~and~~}} a_1=1,a_2=0.$$

\noindent{} From the construction one may check that the sequence
${(T^{h_k})}_{k \in \N^*}$ converges weakly to $\displaystyle
\sum_{k=0}^{+\infty}\frac1{2^k}T^k.$ Thus it is easy to get from
Lemma \ref{beurlingmixing} the following result.

\begin{cor}
The historical Chacon map is spectrally disjoint from any mixing map.
\end{cor}

\begin{rem}Combining the same methods with the continuous version of \linebreak Beurling theorem we can extend Theorem \ref{beurling} to flows. Now, applying the  results of Fr\c{a}czek-Lema\'nczyk \cite{frac-lem1}, \cite{frac-lem2}, we can exhibit examples of maps for which the conditions of the theorem are satisfied.
\end{rem}

\section{Mathew-Nadkarni transformation}
In this section we show that $\alpha$-rigidity alone is not enough to ensure the spectrum is singular. Indeed, we will show that the Mathew-Nadkarni transformation is an example of $\frac12$-rigid transformation whose spectrum has a Lebesgue component.\\

The Mathew-Nadkarni transformation is a $\Z_2$-extension of the
odometer or von Neumann-Kakutani adding machine $\tau$
\footnote{see \cite{Petersen}, \cite{Friedman2} or
\cite{Nadkarni}.}, $\Z_2 \stackrel{\rm {def}}{=} \{0,1\}$ endowed
with its Haar measure $h\stackrel{\rm {def}}{=}
\frac12\delta_0+\frac12\delta_1$. Explicitly, define
$T_\phi~:~[0,1) \times \Z_2 \mapsto [0,1) \times \Z_2$ by
$T_\phi(x,g)=(Tx,\phi(x)+g),$ where $T$ is the
adding machine.\\
$T$ is defined by mapping the interval
$[1-\frac1{2^{n}},1-\frac1{2^{n+1}})$ linearly onto the interval
$[\frac{1}{2^{n+1}},\frac{1}{2^n})$, $n=0,1,2\cdots .$ Note that $T$ is a
rank one transformation which is also easily described using
cutting and stacking. The cocycle $\phi$ is defined (inductively
over all "levels" of the tower associated to $T$ except the last
one) to be $0$ on the intervals
$[1-\frac1{2^{n}},1-\frac1{2^{n}}+\frac1{2^{n+2}})$, and 1 on the
intervals
$[1-\frac1{2^{n}}+\frac1{2^{n+2}},1-\frac1{2^{n+1}})$, $n=0,1,2\cdots .$.\\

\begin{thm} [Mathew-Nadkarni]
The Mathew-Nadkarni transformation has a spectrum consisting of a
Lebesgue component with multiplicity 2, together with discrete
component.
\end{thm}
The proof can be found in {\cite {Mathew-Nadkarni}} or {\cite{Lem}}. We will proof the following theorem.
\begin{thm}
The Mathew-Nadkarni transformation is $\displaystyle
\frac12$-rigid.
\end{thm}

\begin{pf}
Recall that the operator $U_{T_\phi}$ : $L^2([0,1) \times \Z_2,\mu
\otimes h) \rightarrow L^2([0,1) \times \Z_2,\mu \otimes h) $ has
a direct sum decomposition
\[
L^2([0,1) \times \Z_2,\mu \otimes h)=L_0 \oplus L_1,
\]
\noindent where $L_0=\{f\otimes1~:~  f \in L^2([0,1))\}$ and $L_1=\{f
\otimes \chi ~:~ f \in L^2([0,1)) {\rm {~and~}} \chi(g)=(-1)^g,
 {\rm {~for~}} g \in \Z_2\}.$\\
Let $A$ be a Borel set and $\varepsilon \in \Z_2$ then
\[
\1_{A \times \{\varepsilon\}}=f_1+\chi~f_2,
\]
\noindent where $\displaystyle f_1=\frac12 \1_A, f_2=\chi(\varepsilon)f_1$.
Let $\sigma_d$ be the discrete part of the spectral measure
$\sigma_{\1_{A \times \{\varepsilon\}}}$. Then that there exists
a sequence of integers $(n_k)$ such that
\[
\lim_{k \rightarrow \infty}\stackrel{\wedge}{\sigma}_{A \times
\{\varepsilon\}}(n_k)=\lim_{k \rightarrow
\infty}\stackrel{\wedge}{\sigma}_{d} (n_k)=(f_1,f_1)=\frac12 m
\otimes h(A \times \{\varepsilon\}).
\]
Hence, by a standard argument, for any Borel set $B$ in
the $\sigma$-algebra of $[0,1) \times \Z_2$ we have
\[
\liminf_{k \rightarrow \infty}\stackrel{\wedge}{\sigma}_{B}(n_k)
\geq \frac12 m \otimes h(B).
\]
This completes the proof of the theorem.
\end{pf}

\begin{rem}
The Mathew-Nadkarni transformation is not weakly mixing but using
the Ageev's construction \cite {Ageev} one can produce an
$\alpha$-rigid transformations with continuous spectrum and
Lebesgue component of even multiplicity. Let us remark also that
Mathew-Nadkarni's construction contains continuum pairwise
non-isomorphic dynamical systems \cite {Lem}. It follows that one
can produce a continuum of $\alpha$-rigid transformations with
2-fold Lebesgue spectrum.\\
V. V. Ryzhikov told us that Ageev had obtained $\alpha$-rigid maps with Lebesgue spectrum using the examples of Parry-Helson. This result has never been published. Ryzhikov also communicated us that Ageev used the
same assertion as it was
 pointed out to us also by M. Lemanczyk:\\
If $T$ is rigid and the cocyle $\phi : X \longrightarrow \Z_2$ gives Lebesgue spectrum for $T_{\phi}$ (such $\phi$ exists over each rigid maps by Helson-Parry \cite{Helson-Parry}), then $T_{\phi}$ is $\frac12$-rigid.

\end{rem}

\section{Substitution examples}

In this section we review some examples of $\alpha$-rigid transformations whose spectrum has Lebesgue component. In particular we will
give an example with $0<\alpha\ll\frac12$ coming from substitution theory.\\

A vast literature is devoted to substitutions, whose Bible is
\cite{Queffelec2}. They appear as symbolic systems defined on a
finite alphabet $A=\{0,1,\cdots,k-1\}$; a substitution $\xi$ is a
mapping from $A$ to the set $A^*$ of all finite words of A. It
extends naturally into a morphism of $A^*$ by concatenation. We
restrict ourselves to the case when $\xi (0)$ begins with $0$ and
the length of $\xi^n (0)$ tends to infinity with $n$. The infinite
sequence $u$ beginning with $\xi^n (0)$ for all $n$ is then called
a fixed point of $\xi$ and the symbolic system
associated to $u$ is called the dynamical system associated to $\xi$.\\

When $\xi$ is primitive, (i.e., there exists $n$ such that $\alpha$
appears in $\xi^n(\beta)$ for all $\alpha, \beta \in A$), the system
is uniquely ergodic, and we can consider the measure-preserving
system associated to $\xi$, to which we refer in short by ``the
substitution $\xi$''. The composition matrix $M$ of $\xi$ is the matrix whose entries are
$\ell_{ij}=O_i(\xi(j))$, where $i,j \in A$ and $O_i(\xi(j))$ is
the number of $i'$s occurring in $\xi(j)$. If $\xi$ is primitive,
it follows from the Perron-Frobenius theorem that $M$ admits a
strictly positive simple eigenvalue $\theta$, such that $\theta >
|\lambda|,$ for any other eigenvalue $\lambda$ and there exists a
strictly positive eigenvector corresponding to $\theta$. It is
easy to see that (see \cite{Queffelec2}), for any $a \in
A$, the sequence of $k$-dimensional vector
$$(\frac{O_0(\xi^n(a))}{\theta^n},\cdots,\frac{O_k(\xi^n(a))}{\theta^n}),$$
converges to a strictly positive eigenvector $v(a)$ corresponding
to $\theta$.

A classical result on primitive substitution is the Keane-Dekking theorem \cite{Keane}. They
proved that the dynamical system arising from a primitive substitution is not mixing. Their proof, with the above notations, contains the following.\\

\begin{thm}
The dynamical system arising from a primitive substitution $\xi$ is
$r\rho$-rigid. Here $r$ is the maximum of the measure of the
cylinders set $[aa]$, $a \in A$, and $\rho$ is the $\ell_1$-norm
of the strictly positive eigenvector $v(a_r)$ corresponding to the
Perron-Frobenius eigenvalue $\theta$ of $\xi$, $a_r$ is a letter
for which the measure of $[a_ra_r]$ is $r$.
\end{thm}

\noindent Actually it is obvious to compute the constant $r$. In fact, let $\xi_2$ be the substitution defined on the alphabet $A_2=\{(ab),~a,b \in A\}$ in the following way :\\
$${\rm {if}}~~~~ \xi(ab)=\xi(a)\xi(b)=y_0y_1y_2y_3,$$\\
\noindent then we set
$$\xi_2(ab)=(y_0y_1)(y_1y_2).$$
\noindent Now it is easy to compute the normalized positive eigenvector
corresponding to the dominant eigenvalue of the composition matrix
$M_2$ of $\xi_2$.

One of the classical examples of substitutions is the
Rudin-Shapiro substitution. This is
defined on the alphabet $A=\{0,1,2,3\}$ in the following way
\begin{eqnarray*}
\xi(0)=02,~~&&\xi(1)=32,\\
\xi(2)=01,~~&&\xi(3)=31.
\end{eqnarray*}
M. Queff\'elec in \cite{Queffelec1} shows that the continuous part
of the Rudin-Shapiro dynamical system is Lebesgue with
multiplicity 2 and it is easy to prove that the Rudin-Shapiro
dynamical system is $\displaystyle\frac12$-rigid. Another example
of a substitution with Lebesgue spectrum is given by the following
substitution $\xi$ on the alphabet $A=\{0,1,2\}$:
\begin{eqnarray*}
\xi(0)&=&001\\
\xi(1)&=&122\\
\xi(2)&=&210.
\end{eqnarray*}
This substitution has Lebesgue spectrum with multiplicity 2 in the
orthocomplement of eigenfunctions \cite[p.~221]{Queffelec2}. We
point out that one may use a standard computer program to
compute approximatively the constant of $\alpha$-rigidity of $\xi$ which is
approximatively equal to $0.3104979673 \times 10^{-7}.$\\

\begin{thank}
 I would like to express my thanks to Late J. De Sam
Lazaro \footnote{ This paper was prepared under the supervision of
professor Jos\'e De Sam Lazaro two years before his death.} , J-P. Thouvenot and M. Disertori for their considerable help and encouragement.
I'm greatly indebted to V. V. Ryzhikov for the proof of Baxter theorem and for pointing out to me Lemma \ref{beurlingmixing}.
Finally I wish to express my heartfelt gratitude to F. Parreau for his collaboration in proving Theorem 10 and  to M. Lema\'nczyk for  stimulating conversations on the subject and the many suggestions for applications of our main result.

\end{thank}

\bibliographystyle{nyjplain}

\begin{thebibliography}{A-O-W}

\bibitem{elabdal}
E. ~H. ~el ~Abdalaoui, ~A. ~Nogueira, ~T.~de ~la ~Rue, Weak mixing of maps with bounded cutting parameter. {\em  New York J. Math.}, {\bf {11}},  (2005), 81--87 (electronic).

\bibitem{Ageev}
O.~N. ~Ageev,{ Dynamical systems with an even-multiplicity
Lebesgue component in the spectrum}, {\em Math. USSR. Sbornik.}, {\bf
Vol 64 (2)}, (1989), 305--316.

\bibitem{Arnoux}
P.~Arnoux, ~D. ~Ornstein, ~B. ~Weiss, { Cutting and stacking,
interval exchanges and geometric models}, {\em Isr. J. Math.}, {\bf Vol
50 (1-2)}, (1985), 160--168.

\bibitem{Beurling}
A. ~Beurling, { Quasianalyticity and General Distributions},
{\em Stanford Univ. Lect. Notes}, no. 30, {1961.}

\bibitem{Baxter}
 J. ~R.~Baxter, {On the Class of Ergodic
Transformations},  {\em PhD Thesis, University of Toronto}, 1969.

\bibitem{Chacon1}
R.~V. ~Chacon, { Weakly mixing transformations which are not strongly mixing,}
{\em  Proc. Amer. Math. Soc.,} {\bf {22}} (1969), 559--562.

\bibitem{Chacon2}
R.~V. ~Chacon, {Transformations having continuous spectrum,}
{\em J. Math. Mech.,} {\bf  16}  (1966), 399--415.

\bibitem{Keane}
F. ~M. ~Dekking ~and ~M.~Keane, { Mixing properties of
substitutions,} {\em Zeit.~Wahr.}, {\bf 42} (1978), 23--33.

\bibitem{Del}
A. del Junco~ and ~M. ~Lema\'nczyk, {Generic spectral
properties of measure-preserving maps and applications, } {\em Proc.
Am. Math. Soc.}, {\bf 115} (1992), 725--736.

\bibitem{Del-Ru}
A. ~del Junco~ and ~D. ~Rudolph, { On ergodic actions whose
self-joinings are graphs,} {\em Ergodic. Th and Dynam. Sys.}, {\bf 7},
(1987), 531--557.

\bibitem{Deljunco}
A. ~del Junco, A.~Rahe, L.~ Swanson, {Chacon's automorphism
has minimal self-joinings,} {\em J. Anal. Math.},{\bf 37} (1980),
276--284.

\bibitem{frac-lem1}
K.~Fra\c czek, M.~ Lema\'nczyk, {On disjointness properties of some smooth flows,} {\em Fund. Math.}, {\bf 185}  (2005),  no. 2,
117--142.

\bibitem{frac-lem2}
K.~Fra\c czek, M.~ Lema\'nczyk, {A class of special flows over irrational rotations which is disjoint from mixing flows,}
{\em Ergodic Theory Dynam. Systems,} {\bf 24}  (2004),  no. 4, 1083--1095.

\bibitem{Friedman1}
N. Friedman~, { Partial mixing, partial rigidity, and factors,
} {\em Contemp. Math.},{\bf 94,} Collection: Measure and measurable
dynamics (Rochester, NY, 1987), 141--145.

\bibitem{Friedman2}
N. ~A. ~Friedman, { Introduction to Ergodic Theory}, {\em Van
Nostrand Reinhold}, New York, 1970.

\bibitem{Friedman3}
N.~A. Friedman, {Replication and stacking in ergodic
theory, }  {\em Amer. Math. Monthly, } {\bf 99} (1992), 31--34.


\bibitem{guenais}
M.~Guenais, {Singularit\'e des produits de Anzai associés aux
fonctions caractéristiques d'un intervalle.}  {\em Bull. Soc.
Math. France},   {\bf {127}}  (1999),  no. 1, 71--93.

\bibitem{Helson-Parry}
H. ~Helson ~and ~W. ~Parry, {Cocycles and spectra, } {\em Arkiv
Math.},  {\bf 16} (1978), 195--206.

\bibitem{Mathew-Nadkarni}
J. ~Mathew ~and ~M. ~G. ~Nadkarni, {A measure-preserving
transformation whose spectrum has a Lebesgue component of
multiplicity two, } {\em Bull. Lond. Math. Soc.},  {\bf 16} (1984),
402--406.

\bibitem{Nadkarni}
M.~G. ~Nadkarni, {Spectral Theory of Dynamical Systems },
{\em Birkh\"auser Advanced Texts: Basel Textbooks, Birkh\"auser Verlag}, Basel, 1998.

\bibitem{Oseledets}
V.~I.~Oseledets, { On the spectrum of ergodic automorphisms,} {\em Doklady Akad. Nauk SSSR},
{\bf 168, 5} (1966),1009-1011 (in Russian), translated in {\em Soviet Math. Doklady}, {\bf 7} (1966), 776--779.

\bibitem{Parry}
W. ~Parry, {  Topics in Ergodic Theory }, {\em Cambridge University
Press}, 1981.

\bibitem{Petersen}
 K. ~Petersen, {Ergodic Theory }, {\em Cambridge University Press}, 1983.

\bibitem{Prikhod'ko-Ryzhikov}
A. A. Prikhod'ko ~and ~V. ~V. ~Ryzhikov,{
 Disjointness of the convolutions for Chacon's
automorphism, } Dedicated to the memory of Anzelm Iwanik. {\em Colloq.
Math.}, {\bf 84/85} (2000), {\bf part 1}, 67--74.

\bibitem{Kamae}
 T. Kamae, { Spectral properties of automata generating sequences , }
  {\em Unpublished}.

\bibitem{Katok1}
A.~B. Katok, { Constructions in Ergodic Theory, } {\em preprint} .

\bibitem{Katok2}
 A.~B. Katok, { Interval exchange transformations and some special flows are not mixing, }
 {\em Isr. J. Math.}, {\bf Vol 35}, 4, (1980), 301--310.

\bibitem{Katok3}
A.~B.~Katok,~A.~M.~Stepin, {Approximation of ergodic dynamic systems by periodic transformations,} {\em Dokl. Akad. Nauk SSSR,} {\bf 171} (1966), 1268-1271;
{\em Soviet Math. Dokl.,} {\bf 7} (1966), 1638--1641.

\bibitem{King-Thouvenot}
J.~L.~King and ~J-P.~Thouvenot, { A canonical structure theorem
for joining-rank maps, }  {\em J. Analyse Math.},  {\bf 51 } (1988),
182--227.

\bibitem{Klemes-Reinhold}
I. ~Klemes \& ~K. ~Reinhold, { Rank one transformations with
singular spectre type}, {\em Isr. J. Math.},{\bf vol 98}, (1997), 1--14.

\bibitem{Lem}
M. ~Lema{\'n}czyk, { Toeplitz $\Z_2$-extensions, }
 {\em Ann. Inst. Henri Poincar\'e},  {\bf vol 24,} n$^\circ$ 1, (1988), 1--43.

\bibitem{Queffelec1}
M. ~Queff\'elec, { Une nouvelle propri\'et\'e des suites de
Rudin-Shapiro }{\em t Ann. Fourier, Grenoble}, {\bf 37}, 2 (1987), 115--138.

\bibitem{Queffelec2}
M. ~Queff\'elec, {Substitution Dynamical Systems-Spectral
Analysis}, {\em In A. Dold and B. Eckmann, editors, Lecture notes
 in Mathematics}, {\bf {vol 1294}}, Spring-Verlag, 1987.

\bibitem{Rudin}
W. ~Rudin, {Fourier Analysis on Groups }, {\em Interscience
 Tracts in Math.} 12. J. Wiley ans sons, 1962.

\bibitem{Robinson}
E.~A.~,~Jr.~Robinson, {Ergodic measure preserving transformations with arbitrary finite spectral multiplicities,}
{\em Invent. Math.,}  {\bf 72}  (1983),  no. 2, 299--314.

\bibitem{Rokhlin}
V. ~A. ~Rokhlin,  {Selected topics in the metric theory of
dynamical systems, } { \em Uspekhi Mat. Nauk. " New series"}, {\bf 4}
(1949), 57--128 (Russian); {\em  Amer. Math. Soc. Transl.}, 2, {\bf 40}
(1966), 171--240, 1979,  97--122.


\bibitem{Rudolph}
D. ~Rudolph,  { An Example of a measure-preserving map with
minimal self-joining and applications, }
 { \em J. analyse Math.}, {\bf 35},
1979,  97--122.

\bibitem{Ryzhikov}
V. ~V. ~Ryzhikov, {Joinings, Intertwining operators,
 factors and mixing properties of dynamical systems },
 {\em Russian. Acad. Sci. Izv. Math.},  {\bf Vol 26}, (1994), 91--114.

\bibitem{RyzhikovS}
V. ~V. ~Ryzhikov, {Around simple dynamical systems. induced
joinings and multiple mixing}
{\em J. Dynam. Control Systems,} {\bf {3}}  (1997),  no. 1, 111--127.

\bibitem{RyzhikovP}
V. ~V. ~Ryzhikov, {Stochastic intertwinings and multiple mixing of dynamical systems,}
{\em J. Dynam. Control Systems,} {\bf {2}}  (1996),  no. 1, 1--19.

\bibitem{Ryzhikovemail}
V. ~V. ~Ryzhikov, Prive communication.

\bibitem{Stepin}
A.~M. Stepin, {Spectral properties of generic dynamical
system, } {\em Math. U.S.S.R. Izv.}, {\bf 50} (1986), 159--192.

\bibitem{Veech}
W. ~A. ~Veech, {The metric theory of interval exchange transformations. I Generic spectral
properties,} {\em Amer. J. Math.}, {\bf 106} (1984), 1331--1359.

\bibitem{Veech2}
W. ~A. ~Veech, { A criterion for a process to be prime, }
 {\em Monatsh. Math.},  {\bf 94} (1982), 335--341.

\end{thebibliography}

\centerline{(received version 10.10.2007)}
\centerline{(revised version 04.01.2009)}

Author's address:

E.~H. ~ El Abdalaoui

Department of Matehmatics,
LMRS UMR 60 85 CNRS,

University of Rouen,

Avenue de l'Universit\'e, BP.12,
76801 Saint Etienne du Rouvray, France

E-mail:elhoucein.elabdalaoui@@univ-rouen.fr

\end{document}